\documentclass{amsart}[12pt]
\usepackage{amsmath,amsfonts,amssymb,latexsym,amsthm}

\usepackage{relsize,etoolbox}
\AtBeginEnvironment{quote}{\smaller}

\def\al{\alpha}

\def\<{\langle}
\def\>{\rangle}

\def\ts{\hskip.015cm}

\def\0{{\mathbf 0}}

\def\nin{\noindent}
\def\sm{\smallskip}
\def\ms{\medskip}
\def\.{\hskip.06cm}
\def\ts{\hskip.03cm}

\begin{document}
\title{History of Catalan numbers}

\author[Igor~Pak]{ \ Igor~Pak$^\star$}


\begin{abstract}
We give a brief history of Catalan numbers, from their first discovery
in the 18th century to modern times. This note will appear as an appendix
in Richard Stanley's forthcoming book~\cite{Sta2}.
\end{abstract}

\thanks{\thinspace ${\hspace{-.45ex}}^\star$Department of Mathematics,
UCLA, Los Angeles, CA, 90095.
\hskip.06cm
Email:
\hskip.06cm
\texttt{pak@math.ucla.edu}}

\maketitle

\vskip.9cm

\section*{Introduction}\label{intro}

\noindent
In the modern mathematical literature, Catalan numbers are wonderfully ubiquitous.
Although they appear in a variety of disguises, we are so used to having them
around, it is perhaps hard to imagine a time when they were either unknown
or known but obscure and underappreciated.  It may then come as a surprise
 that Catalan numbers have a rich history of multiple rediscoveries until
 relatively recently.  Below we review over 200 years of history, from their
 first discovery to modern times.

 \smallskip

We break the history into short intervals of mathematical activity,
each covered in a different section.  We spend most of our effort on
the early history, but do bring it to relatively recent developments.  
We should warn the reader that although this work is in the History of 
Mathematics, we are not a mathematical historian.  Rather, this work is 
more of a historical survey with some added speculations based on our
extensive reading of the even more extensive literature.  Due to
the space limitations, this survey is
very much incomplete, as we tend to emphasize first discoveries
and papers of influence rather than describe subsequent developments.

This paper in part is based on our earlier investigation reported
in~\cite{Pak-blog}.  Many primary sources are assembled on the
\emph{Catalan Numbers website}~\cite{Pak}, including scans of the
 original works and their English translations.

\medskip

 \section{Ming Antu}

\nin
Ming Antu, sometimes written Ming'antu (c.1692--c.1763)
was a Chinese scientist and mathematician, native of Inner Mongolia.
In 1730s, he wrote a book \emph{Quick Methods for Accurate Values
of Circle Segments} which included a number of
trigonometric identities and power series, some involving
Catalan numbers:
$$
\sin(2\alpha) \, = \, 2 \sin \al \, - \,\sum_{n=1}^\infty
\.  \frac{C_{n-1}}{4^{n-1}} \. \sin^{2n+1}\alpha \, = \,
2\sin\alpha \ts - \ts \sin^3\alpha \ts - \ts \frac14\sin^5\alpha
\ts - \ts \frac18 \sin^7\alpha \ts - \ts \ldots
$$
The integrality of Catalan numbers played no role in this work.

Ming Antu's book was published only in 1839, and the connection to
Catalan numbers was observed by Luo Jianjin in 1988. We refer
to~\cite{Lar1,Luo} for more on this work and further references.

\medskip

\section{Euler and Goldbach}

\nin
In 1751, Leonhard Euler (1707--1783) introduced and found a closed formula
for what we now call the \emph{Catalan numbers}.  The proof of this result
had eluded him, until he was assisted by Christian Goldbach (1690--1764),
and more substantially by Johann Segner.  By~1759,
a complete proof was obtained.  This and the next section tell
the story of how this happened.

On September~4, 1751, Euler wrote a letter to Goldbach which
among other things included his discovery of Catalan numbers.\footnote{The
letters between Euler and Goldbach were published in~1845 by P.~H.~Fuss,
the son of Nicolas Fuss.
The letters were partly excised, but the originals also survived. Note that
the first letter shows labeled quadrilateral and pentagon figures missing
in the published version, see~\cite{Pak}.}
Euler was in Berlin (Prussia) at that time, while his friend and
former mentor Goldbach was in St.~Petersburg (Imperial Russia).
They first met when Euler arrived to St.~Petersburg back in 1727
as a young man, and started a lifelong friendship, with 196
letters between them~\cite{Var}.

Euler defines Catalan numbers $C_n$ as the number of triangulations of
$(n+2)$-gon, and gives the values of $C_n$ for $n\le 8$ (evidently,
computed by hand).  All these values, including $C_8=1430$ are
correct. Euler then observes that successive ratios have a pattern
and guesses the following formula for Catalan numbers:
\begin{equation}\label{eq:app-product}
C_{n-2} \, = \, \frac{2\cdot 6 \cdot 10 \cdots (4n-10)}
{2 \cdot 3 \cdot 4 \cdots (n-1)} \,.
\end{equation}
For example, $C_3=(2\cdot 6 \cdot 10)/(2\cdot 3 \cdot 4) =5$.
He concludes with the formula for the Catalan numbers g.f.~:
\begin{equation}\label{eq:app-gf}
A(x) \. = \. 1 \ts + \ts 2x  \ts + \ts 5x^2 \ts + \ts 14x^3  \ts + \ts 42x^4\ts + \ts 132x^5  \ts +\ts \ldots
\. = \. \frac{1-2x-\sqrt{1-4x}}{2x^2}\..
\end{equation}

In his reply to Euler dated October 16, 1751, Goldbach notes that
the g.f.~$A(x)$ satisfies quadratic equation
\begin{equation}\label{eq:app-quadratic}
1\ts + \ts x A(x) \. = \. A(x)^{\frac12}\..
\end{equation}
He then suggests that this equation can be used to derive Catalan
numbers via an infinite family of equations on its coefficients.

Euler writes back to Goldbach on December 4, 1751.  There, he explains
how one can obtain the product formula~\eqref{eq:app-product} from the
binomial formula:
\begin{equation}\label{eq:app-binomial}
\sqrt{1-4x}\. = \. 1 \ts - \ts \frac12 \. 4x \ts - \ts \frac{1\cdot 1}{2\cdot 4} \. 4^2x^2
 \ts - \ts \frac{1\cdot 1\cdot 3}{2\cdot 4\cdot 6} \. 4^3x^3 \ts -
 \ts \frac{1\cdot 1\cdot 3  \cdot 5}{2\cdot 4\cdot 6\cdot 8} \. 4^4x^4 \ts - \ts \ldots
\end{equation}
From the context of the letter, it seems that Euler knew the exact form
of~\eqref{eq:app-binomial} before his investigation of the Catalan numbers.
Thus, once he found the product formula~\eqref{eq:app-product},
he was able to derive~\eqref{eq:app-gf} rather than simply guess it.
We believe that Euler did not include his derivation in the first letter
as the g.f.~formula appears at the very end of it,
but once Goldbach became interested he patiently
explained all the steps, as well as other similar formulas.

\medskip

\section{Euler and Segner}

\nin
Johann Andreas von Segner (1704--1777) was another frequent correspondent
of Euler.  Although Segner was older, Euler rose to prominence faster, and
in 1755 crucially helped him to obtain a position at University of Halle,
see~\cite[p.~40]{Cal}.  It seems, there was a bit of competitive tension
between them, which adversely affected the story.

In the late 1750s, Euler suggested to Segner the problem of counting
the number of triangulations of an $n$-gon.  We speculate\footnote{Euler's letters
to Segner did not survive, as Segner directed all his archive to be burned
posthumously~\cite[p.~153]{Fell}. Segner's letters to Euler did
survive in St.~Petersburg, but have yet to be digitized (there are 159 letters
between 1741 and 1771).}
 based on Segner's later work, that Euler told him only values up to~$C_7$, and neither
the product formula~\eqref{eq:app-product} nor the g.f.~\eqref{eq:app-gf}.

Segner accepted the challenge and in 1758 wrote a paper~\cite{Segner},
whose main result is a recurrence relation which he finds and proves
combinatorially:
\begin{equation}\label{eq:app-segner}
C_{n+1}\. = \. C_0C_n \ts + \ts  C_1C_{n-1} \ts + \ts C_2C_{n-1} \ts +
\ts \ldots  \ts + \ts C_nC_0\ts.
\end{equation}
He then uses the formula to compute the values of~$C_n$, $n\le 18$,
but makes an arithmetic mistake in computing $C_{13} =742,900$, which then
invalidates all larger values.

Euler must have realized that equation~\eqref{eq:app-segner} is the last
missing piece necessary to prove~\eqref{eq:app-quadratic}.
He arranged for Segner's paper to be published in the journal of
St.~Petersburg Academy of Sciences, but with his own\footnote{The
article is unsigned, but the authorship by Euler is both evident
and reported by numerous sources.}    \emph{Summary}~\cite{Euler}.
In it, he states~\eqref{eq:app-product}, gives Segner a lavish
compliment, then points out his numerical mistake, and correctly
computes all $C_n$ for $n\le 23$.  It seems unlikely
that given the simple product formula, Segner would have computed
$C_{13}$ incorrectly, so we assume that Euler shared it only with his close
friend Goldbach, and kept Segner in the dark until after the publication.

In summary, a combination of results of Euler and Segner, combined
with Goldbach's observation gives a complete proof of the product
formula~\eqref{eq:app-product}.  Unfortunately, it took about
80~years until the first complete proof was published.

\medskip

\section{Kotelnikow and Fuss}

\nin
Sem\"{e}n Kirillovich Kotelnikow (1723--1806) was a Russian
mathematician of humble origin who lived in St.~Petersburg
all his life.  In 1766, soon after the Raid on Berlin, Euler
returned to St.~Petersburg.  Same year, Kotelnikow writes a
paper~\cite{Kot} elaborating on Catalan numbers.  Although he
claimed to have another way to verify~\eqref{eq:app-product},
Larcombe notes that he ``does little more than play
around with the formula''~\cite{Lar2}.

Nicolas Fuss (1755--1826) was a Swiss born mathematician
who moved to St.~Petersburg to become Euler's assistant
in~1773.  He married Euler's granddaughter, became a well
known mathematician in his own right, and remained
in Russia until death.  In his 1795 paper~\cite{Fuss},
in response to Pfaff's question on the number of subdivisions
of an $n$-gon into $k$-gons, Fuss introduced what is
now know as the \emph{Fuss--Catalan numbers} (see Exercise~A14), and
gave a generalization of Segner's formula~\eqref{eq:app-segner},
but not the product formula.

\medskip

 \section{The French school, 1838--1843}

\nin
In 1836, a young French mathematician
Joseph Liouville (1809--1882) founded the
\emph{Journal de Math\'{e}matiques Pures et Appliqu\'{e}es}.
He was in the center of mathematical life in Paris
and maintained a large mailing list, which proved
critical in this story.

In 1838, a Jewish French mathematician and mathematical historian
Olry Terquem (1782--1862) asked Liouville if he knows a simple
way to derive Euler's formula~\eqref{eq:app-product} from Segner's
recurrence~\eqref{eq:app-segner}.   Liouville in turn communicated
this problem to ``various geometers''.   What followed is
a remarkable sequence of papers giving foundation to 
``Catalan studies''.

First, Gabriel Lam\'{e} (1795--1870) wrote a letter to Liouville
outlining the solution, a letter Liouville promptly published in
the \emph{Journal} and further popularized~\cite{Lame}.  Lam\'{e}'s solution
was to use an elegant double counting argument.  Let's count the
number $A_n$ of triangulations of an $(n+2)$-gon with one of its
$(n-1)$ diagonals oriented.  On the one hand,
$A_n=2(n-1) C_n$.  On the other hand, by summing over all
possible directed diagonals we have
$$
A_n \, = \, n \ts \bigl(C_1C_{n-1} + C_2C_{n-2} + \ldots + C_nC_{1}\bigr)\ts.
$$
Combining these two formulas with~\eqref{eq:app-segner} easily
implies~\eqref{eq:app-product}.

In 1738, Belgium born mathematician Eug\`{e}ne Charles Catalan  (1814--1894)
was a student of Liouville at \'{E}cole Polytechnique.
Inspired by the work of Lam\'{e}, he became interested in the problem.
He was the first to obtain what are now standard formulas
$$
C_n \. = \. \frac{2n!}{n! \ts (n+1)!} \. = \.
\binom{2n}{n} -  \binom{2n}{n-1}\ts.
$$
He then studied the problem of computing the number of different
(non-associative) products of $n$ variables, equivalent to counting
the number of bracket sequences~\cite{Cat1}.

Olinde Rodrigues (1795--1851) was a descendant of a large Sephardic
Jewish family from Bordeaux.  He received his doctorate in
mathematics and had a career as a banker in Paris, but continued his
mathematical interests.  In the same volume of the \emph{Journal},
he published two back-to-back short notes giving a more direct
double counting proof of~\eqref{eq:app-product}.  The first
note~\cite{Rod1} gives a variation on Lam\'{e}'s argument, a
beautiful idea often regarded a folklore.  Roughly,
he counts in two ways the number $B_n$ of triangulations of $(n+2)$-gon where
either an edge or a diagonal is oriented.  On the one hand,
$B_n = 2(2n+1) C_n$.  On the other hand, $B_n= (n+2)C_{n+1}$,
since triangulations with an oriented diagonal are in bijection
with triangulations of $(n+3)$-gon obtained by inserting a triangle
in place of a diagonal, and such edge can be any edge except the first one.
We omit the details.

In~\cite{Rod2}, Rodrigues gives a related, but even simpler argument for counting
the bracket sequences.  Denote by $P_n$ the number of bracket sequences
of labeled terms $x_1,\ldots,x_n$, e.g. $x_2(x_1x_3)$.  Then, on the
one hand $P_n=n! \ts C_n$.  On the other hand, $P_{n+1} = (4n-2)P_n$ since
variable $x_{n+1}$ can be inserted into every bracket sequence in exactly
$(4n-2)$ ways.  To see this, place a bracket around every variable
and the whole product, e.g. $((x_2)((x_1)(x_3)))$.
Now observe that a new variable is inserted immediately
to the left of any of the $(2n-1)$ left brackets, or immediately
to the right of any of the $(2n-1)$ right brackets,
 e.g.
 $$
 ((x_2)((x_1\mathbf{\underline{)}}(x_3))) \. \to \. ((x_2)(((x_1)(x_4))(x_3)))\ts.
 $$

In 1839, clearly unaware of Euler's letters, a senior French
mathematician Jacques Binet (1786--1856) wrote a paper~\cite{Binet}
with a complete g.f.~proof of~\eqref{eq:app-product}.
Across the border in Germany, Johann August Grunert (1797--1872)
became interested in the work of the Frenchmen on the one hand,
and of Fuss on the other.  A former student of Pfaff, in a 1841
paper~\cite{Gru} he found a product formula for the Fuss--Catalan
numbers.  He employed g.f.'s to reduce the problem to
$$
Z(x)^m \, = \, \frac{Z(x)-1}{x}\.,
$$
but seemed unable to finish the proof.  The complete proof was given in 1843
by Liouville~\cite{Liu} using the
\emph{Lagrange inversion} (see~\cite{LW1}).\footnote{Liouville was clearly
unaware of Fuss's paper until~1843, yet gives no credit to Grunert (cf.~\cite{LW}).}

After a few more mostly analytic papers inspired by the problem, the
attention of the French school turned elsewhere.  Catalan however returned
to the problem on several occasions throughout his career. Even fifty years later,
in 1878, he published \emph{Sur les nombres de Segner}~\cite{Cat2} on divisibility
of the Catalan numbers.

\medskip

\section{The British school, 1857--1891}

\nin
Rev.~Thomas Kirkman (1806--1895) was a British ordained minister who
had a strong interest in mathematics.  In 1857, unaware of the previous work
but with Cayley's support, he published a lengthy treatise~\cite{Kirkman}.
There, he introduced the \emph{Kirkman--Cayley numbers}, defined as the number
of ways to divide an $n$-gon with $k$ non-intersecting diagonals
(see~A41), and states a general product formula, which he proves
in a few special cases.

In 1857, Arthur Cayley (1821--1895) was a lawyer in London and extremely
prolific mathematically.  He was interested in a related counting of
\emph{plane trees}, and in 1859 published a short note~\cite{Cayley1},
where (among other things) he gave a conventional g.f.\ proof of
that the number of plane trees is the Catalan number.  Like Kirkman,
he was evidently unaware of the previous work.  His final formula is
$$
C_{m-1} \, = \, \frac{1\cdot 3 \cdot 5 \cdots (2m-3)}{1\cdot 2 \cdot 3 \cdots m} \. 2^{m-1}\.,
$$
which he called ``a remarkably simple form''.  Curiously, in the same paper
he also discovered and computed the g.f.\ for the \emph{ordered Bell numbers}.

In 1860, Cayley won a professorship in Cambridge and soon became
a central figure in British mathematics.  Few years later,
Henry Martyn Taylor (1842--1927) became a student in Cambridge,
where he remained much of his life, working in geometry, mathematical
education and politics.\footnote{Despite his blindness, Taylor
was elected a Mayor of Cambridge in 1900; see {\tt http://tinyurl.com/md99y9b}.}
In 1882, he and R.~C.~Rowe published a paper~\cite{TR} where they carefully
examined the literature of the French school, but gave an erroneous
description of the Euler--Segner story.
They evidently missed later papers by Grunert and Liouville, and
computed the Fuss--Catalan numbers along similar lines, using the
Lagrange inversion, see~\cite{LW}.

Cayley continued exploring g.f.~methods for a variety of enumerative
problems, among over 900 papers he wrote.  In 1891, Cayley used g.f.~tour
de force to completely resolve Kirkman's problem~\cite{Cayley2}.

\medskip

\section{The ballot problem}

\nin
The \emph{ballot numbers} were first defined by Catalan in~1839, disguised as
the number of certain triangulations~\cite{Cat3}.  We believe this contribution
was largely forgotten since Catalan gave a formula for the ballot numbers in
terms of the Catalan numbers, but neither gave a closed formula nor even a table
of the first few values.

The ballot sequences were first introduced in 1878 by
William Allen Whitworth (1840--1905), a British mathematician, a priest and a
fellow at Cambridge.  He resolved the problem completely by an elegant
counting argument and found interesting combinatorial applications~\cite{Whi}.
Despite both geographical and mathematical proximity to Cayley, he did not
notice that the numbers $1,2,5,14,\ldots$ he computed are the Catalan numbers.

In modern terminology, the
\emph{ballot problem} was introduced in 1887 by Joseph Bertrand (1822--1900).
In a half a page note~\cite{Ber}, he defined the probability $P_{m,n}$ that 
in an election with $(m+n)$ voters, a candidate who gathered $m$ votes is 
always leading the candidate with $n$ votes.  Bertrand announced that 
one can use induction to show that 
$$
P_{m,n} \. = \. \frac{m-n}{m+n}\..
$$
He famously concludes by saying:
\sm

\begin{quote}
Il semble vraisemblable qu'un r\'{e}sultat aussi simple pourrait
se d\'{e}montrer d'une mani\`{e}re plus directe.
(It seems probable that such a simple result could
be shown in a more direct way.)
\end{quote}
\sm

\nin
Within months, Bertrand's prot\'{e}g\'{e} Joseph-\'{E}mile Barbier (1839--1889)
announced a generalization of $P_{m,n}$ to larger proportions of leading
votes~\cite{Bar}, a probabilistic version of the formula for the
Fuss--Catalan numbers.

D\'{e}sir\'{e} Andr\'{e} (1840--1918) was a well known
French combinatorialist, a former student of Bertrand.  He published an
elegant proof of Bertrand's theorem~\cite{Andre} in the same 1887 volume
of \emph{Comptes Rendu} as Bertrand and Barbier.  Although he is often
credited with the \emph{reflection method}, this attribution is
incorrect,\footnote{Marc Renault was first to note the mistake~\cite{Ren};
he traces the confusion to \emph{Stochastic Processes} by J.L.~Doob (1953)
and \emph{An Introduction to Probability Theory} by W.~Feller (2nd ed., 1957).
In a footnote, Feller writes:  ``The reflection principle is used
frequently in various disguises, but without the geometrical
interpretation it appears as an ingenious but incomprehensible trick. The
probabilistic literature attributes it to D. Andr\'{e} (1887). It appears in connection with
the difference equations for random walks.  These are related to some
partial differential equations where the reflection principle is a familiar tool called
\emph{method of images}. It is generally attributed to Maxwell and Lord Kelvin.''}
as Andre's proof was essentially equivalent to that by
Whitworth (cf.~\cite{Ren}).
The reflection principle is in fact due to Dmitry Semionovitch Mirimanoff
(1861--1945), a Russian-Swiss mathematician who discovered it
in 1923, in a short note~\cite{Mer}.  Andre's original proof was
clarified and extended to Barbier's proposed generalization in~\cite{DM}.

We should mention here Tak\'{a}cs's thorough treatment of
the history of the ballot problem in probabilistic context~\cite{Tak},
Humphreys's historical survey on general lattice paths and the reflection
principle~\cite{Hum}, and Bru's historical explanation on how
Bertrand learned it from Amp\`{e}re~\cite{Bru}.

\medskip

\section{Later years}

\nin
Despite a large literature, for decades the Catalan numbers
remained largely unknown and unnamed in contrast with other
celebrated sequences, such as the Fibonacci and Bernoulli numbers.
Yet the number of publications on Catalan numbers was
growing rapidly, so much that we cannot attempt to cover
even a fraction of them.  Here are few major appearances.

The first monograph citing the Catalan numbers is the
\emph{Th\'{e}orie des Nombres} by \'{E}douard Lucas (1842--1891),
published in Paris in 1891.  Despite the title, the book had
a large combinatorial content, including Rodrigues's
proof of Catalan's combinatorial interpretation in terms of
brackets and products~\cite[p.~68]{Lucas}.

The next notable monograph is the \emph{Lehrbuch der Combinatorik}
by Eugen Netto (1848--1919), published in Leipzig in~1901.   This is one of
the first Combinatorics monographs; it includes a lot of
material on permutations and combinations, and it devotes
several sections to the papers by Catalan, Rodrigues and some
related work by Schr\"oder~\cite[p.~193-202]{Netto}.

We should mention that these monographs were more the exception
than the rule.  Many classical combinatorial books do not
mention  Catalan numbers at all, most notably
Percy MacMahon, \emph{Combinatorial Analysis} (1915/6),
John Riordan, \emph{An Introduction to Combinatorial Analysis} (1958),
and H.~J.~Ryser, \emph{Combinatorial Mathematics} (1963).

A crucial contribution was made by William G.~Brown in~1965,
when he recognized and collected a large number of references
on the ``Euler--Segner problem'' and the ``Pfaff--Fuss problem'',
as he called them~\cite{Brown}.\footnote{Unfortunately, Brown
used somewhat confusing notation $D_{0,m}^{(3)}$ to denote the Catalan numbers.}
From this point on, hundreds of papers involving Catalan numbers have been
published and all standard textbooks started to include it
(see e.g.~\cite[$\S$3.2]{Hall} and~\cite[$\S$3]{vl} for an early adoption).
Some 465 references were assembled  by Henry~W.~Gould in a remarkable
bibliography~\cite{Gould}, which first appeared in~1971
and revised in~2007.   At about that time, various lists of
combinatorial interpretations of Catalan numbers started to
appear, see e.g.~\cite{Gro}, Appendix~1 and a very different
kind of list in~\cite[p.~263]{TF}.

\medskip

\section{The name}

\nin
Because of their chaotic history, the Catalan numbers have received this name
relatively recently.  In the old literature, they were sometimes called
the \emph{Segner numbers} or the \emph{Euler--Segner numbers}, which is
historically accurate as their articles were the first published
work on the subject.  Perhaps surprisingly, we are able to tell exactly
who named them \emph{Catalan numbers} and when.\footnote{This provides
yet another example of the so called \emph{Stigler's Law of Eponymy},
that no result is named after its original discoverer.} Our investigation
of this \emph{eponymy} is informal, but we hope convincing (cf.~\cite{Yc}
for an example of a proper historical investigation).

First, let us discard two popular theories: that the name was introduced
by Netto in~\cite{Netto} or by Bell in~\cite{Bell}.  Upon careful study
of the text, it is clear that Netto did single out Catalan's work and
in general was not particularly careful with the 
references,\footnote{In fairness to Netto, this was standard at the time.}
but never specifically mentioned ``Catalan Zahlen''.
Similarly, Eric Temple Bell (1883--1960) was a well known mathematical
historian, and referred to ``Catalan's numbers'' only in the context
of Catalan's work.  In a footnote, he in fact referred to them as the
``Euler--Segner sequence'' and clarifies the history of the problem.

Our investigation shows that the credit for naming Catalan numbers is
due to an American combinatorialist John Riordan (1903--1988).\footnote{Henry
Gould's response seems to supports this conclusion: \. {\tt http://tinyurl.com/mpyebyw}.}
He tried this three times.  The first two:
\emph{Math Reviews} MR0024411 (1948) and MR0164902 (1964) went unnoticed.
Even Marshall Hall's influential~1967 monograph~\cite{Hall} does not
have the name.  But in 1968, when Riordan used ``Catalan numbers'' in the
monograph~\cite{Rio}, he clearly struck a chord.

Although Riordan's book is now viewed as somewhat disorganized and
unnecessarily simplistic,\footnote{Riordan famously writes in the
introduction to~\cite{Rio}: ``Combinatorialists
use recurrence, generating functions, and such transformations as the
Vandermonde convolution; others, to my horror, use contour integrals,
differential equations, and other resources of mathematical analysis.''}
back in the day it was quite popular.  It was lauded as ``excellent and stimulating''
in P.R.~Stein's review, which continued to say ``Combinatorial identities is,
in fact, a book that must be read, from cover to cover, and several times.''
We are guessing it had a great influence on the field and cemented the
terminology and some notation.

We have further evidence of Riordan's authorship of a different nature:
the \emph{Ngram chart} for the words ``Catalan numbers'',
which searches Google Books content for this sentence over
the years.\footnote{This chart is available here: \.
{\tt http://tinyurl.com/k9nvf28}\ts.}  The search clearly shows that
the name spread after~1968 and that Riordan's monograph was the first book
which contained it.

There were three more events on the way to the wide adoption of the name
``Catalan numbers''.  In 1971, Henry Gould used it in the early version
of~\cite{Gould}, a bibliography of Catalan numbers with 243 entries.
Then, in 1973, Neil~Sloane gave this name to the sequence entry in~\cite{Sloane},
with Riordan's monograph and Gould's bibliography being two of the only
five references.  But the real popularity (as supported by the Ngram chart),
was achieved after the Martin Gardner's \emph{Scientific American} column
in June~1976, popularizing the subject:

\sm

\begin{quote} [Catalan numbers]
have [..] delightful propensity for popping up unexpectedly, particularly in
combinatorial problems.  Indeed, the Catalan sequence is probably the most frequently
encountered sequence that is still obscure enough to cause mathematicians lacking
access to N.J.A.\ Sloane's  \emph{Handbook of Integer Sequences} to expend inordinate
amounts of energy re-discovering formulas that were worked out long ago.~\cite{Gar}
\end{quote}


\medskip

 \section{The importance}

\nin
 Now that we have covered over 200 years of the history of Catalan numbers,
 it is worth pondering if the subject matter is worth the effort.  Here we defer
 to others, and include the following helpful quotes, in no particular order.

\medskip

\nin
Manuel Kauers and Peter Paule write in their 2011 monograph:

\sm

\begin{quote}
It is not exaggerated to say that the Catalan numbers are the most
prominent sequence in combinatorics.~\cite{KP}
\end{quote}

\ms

\nin
Peter Cameron, in his 2013 lecture notes elaborates:
\sm

\begin{quote}
The Catalan numbers are one of the most important sequences of
combinatorial numbers, with a large range of occurrences in apparently
different counting problems.~\cite{Cam}
\end{quote}

\ms

\nin
Martin Aigner gives a slightly different emphasis:

\sm

\begin{quote}
The Catalan numbers are, next to the binomial coefficients,
the best studied of all combinatorial counting numbers.~\cite{Aigner}
\end{quote}

\ms

\nin
This flattering comparison was also mentioned by Neil Sloane and Simon Plouffe in
1995, in the final print edition of~EIS:
\sm

\begin{quote}
Catalan numbers are probably the most frequently occurring combinatorial numbers after
the binomial coefficients.~\cite{SP}
\end{quote}

\ms

\nin
Doron Zeilberger attempts to explain this in his \emph{Opinion~49}~:

\sm

\begin{quote}
Mathematicians are often amazed that certain mathematical objects (numbers, sequences, etc.) show up so often. For example, in enumerative combinatorics we encounter the Fibonacci and Catalan sequences in many problems that seem to have nothing to do with each other.  [..] The answer, once again, is our human predilection for triviality. [..] The Catalan sequence is the simplest sequence whose generating function is a (genuine) algebraic formal power series.~\cite{Zei}
\end{quote}

\ms

\nin
In a popular article, Jon McCammond reveals a new side of Catalan numbers:

\sm

\begin{quote}
The Catalan numbers are a favorite pastime of many amateur (and professional)
mathematicians.~\cite{McC}
\end{quote}

\ms

\nin
Christian Aebi and Grant Cairns give both a praise and a diagnosis:

\sm

\begin{quote}
Catalan numbers are the subject of [much] interest (sometimes known
as \emph{Catalan disease}).~\cite{AC}
\end{quote}

\ms

\nin
Thomas Koshy, in his introductory book on Catalan numbers, gives
them literally a heavenly praise:

\sm

\begin{quote}
Catalan numbers are even more fascinating [than the \emph{Fibonacci numbers}].
Like the North Star in the evening sky,
they are a beautiful and bright light in the mathematical heavens. They continue to
provide a fertile ground for number theorists, especially, Catalan enthusiasts
and computer scientists.~\cite{Koshy}
\end{quote}

\ms

\nin
In conclusion, let us mention the following answer which Richard Stanley gave
in a 2008 interview on how the Catalan numbers exercise came about:

\sm

\begin{quote}
I'd have to say my favorite number sequence is the Catalan numbers. [..]
Catalan numbers just come up so many times. It was well-known before me that they had
many different combinatorial interpretations. [..]
When I started teaching enumerative combinatorics, of course I did the Catalan numbers.
When I started doing these very basic interpretations -- any enumerative course would have
some of this -- I just liked collecting more and more of them and I decided to be systematic.
Before, it was just a typed up list. When I wrote the book, I threw everything I knew in the
book. Then I continued from there with a website, adding more and more problems.~\cite{Kim}
\end{quote}

\vskip.7cm

\noindent
\textbf{Acknowledgments:} \.  We are very grateful to Richard Stanley for
the suggestion to write this note.  We would also like to thank Maria Rybakova
for her help with translation of~\cite{Euler}, to Xavier Viennot for showing
us the scans of Euler's original letter to Goldbach, and to Peter Larcombe 
for sending me his paper and help with other references.  The author was 
partially supported by the~NSF.

\vskip1.2cm


\end{document}